\documentclass[12pt,german]{amsart}

\usepackage{hyperref}

\usepackage{amsmath}
\usepackage{amssymb}
\usepackage{amsthm}
\usepackage{graphicx}
\usepackage{epsfig}

\usepackage{epstopdf}
\usepackage{centernot}

\theoremstyle{definition}

\theoremstyle{remark}

\newcommand{\CC}{{\mathbb C}}

\newcommand{\ff}{\varphi}

\newcommand{\cA}{\mathcal{A}}

\title[]{Some thoughts on Wilhelm von Waldenfels and on universal second order constructions}

\author{Roland Speicher}
\address{Saarland University, Department of Mathematics, D-66123 Saarbr\"ucken, Germany}
\email{speicher@math.uni-sb.de}

\date{\today}

\thanks{Thanks to Wilhelm for showing me that mathematics is hard, but still fun}


\begin{document}

\maketitle

\section{A Tribute to Wilhelm von Waldenfels}

Wilhelm was a master of doing concrete, hands-on calculations and extracting  abstract concepts from such calculations, leading to a deeper understanding of what is going on; as for example in \cite{GivW,GlvW,vWa}. Much of my own work -- not only during my Heidelberg times as a PhD student of Wilhelm's, but more or less during my whole career --  was inspired by this role model.

One example of this spirit were my investigations on the question: how many notions of non-commutative independence are there? Or, from my perspective: how special is free independence? This was something which was started in Heidelberg, mainly in discussions with Michael Sch\"urmann and later resulted in the classification of non-commutative independences in works of Ben Gorbal and Sch\"urmann \cite{BGS}, Muraki \cite{Mur}, and myself \cite{Spe}. It might be fair to say that my contribution is more in Wilhelm's spirit of concrete calculations, whereas Michael took over more of the abstract conceptual approach.

Since those beginnings, when non-commutative life was nice and easy, there have emerged many new versions and generalizations of notions of non-commutative independence, most notably the theory of bi-freeness \cite{Voi} and the one of higher order freeness \cite{MSp}. In such theories, one is not just looking at an algebra equipped with a state (or linear functional), but there is more structure around -- like distinguished subalgebras in the case of bi-freeness or an additional bilinear functional in the case of second order freeness -- and the notion of independence should respect this additional structure. We have in bi-freeness or second order freeness some basic examples for such independences, but again we would like to see how many more are there. Again I would prefer if there are not so many, so that the theories we are working on can be sold as quite unique theories and not just as one possibility in a multiverse of many. 

Recently there has been quite some progress -- in particular by Var\u so \cite{Var} and by Gerhold, Hasebe, and Ulrich \cite{GHU} -- on those issues, most notably for bi-freeness. As to be expected the situation is much more complicated than in the ``classical'' case and there seem to be possibilities at least for non-trivial deformations of the known examples. Maybe there is still hope that modulo such deformations (which might be excluded by some more or less canonical extra normalization conditions) bi-freeness and its few relatives are unique theories. 

From time to time I also try applying my hands-on approach to these questions, in particular, for the second order situation. There are at least some simple observations, which are still far away from a complete classification of such possibilities, but which might be a first step. I think it is in the spirit of Wilhelm if in the following I present  some of these calculations and hope that others might see more in them. 

Let me close this introductory section with another important lesson which I learned from Wilhelm, namely that thinking, in particular on mathematics, is hard and it does not become easier when you get older. At the times in Heidelberg I could not really appreciate this insight and thought that age and experience would help to understand things better and more easily. As I now experience myself this is unfortunately not the case; in order to understand and see some structure, I still have to calculate many examples and hope to get some more insight by discussing also preliminary and inconclusive results with others. Just sitting together with others and talking about your problems often results in unexpected progress. That's what we did in Heidelberg back in the good old days \dots

\quad

\includegraphics[width=4.in]{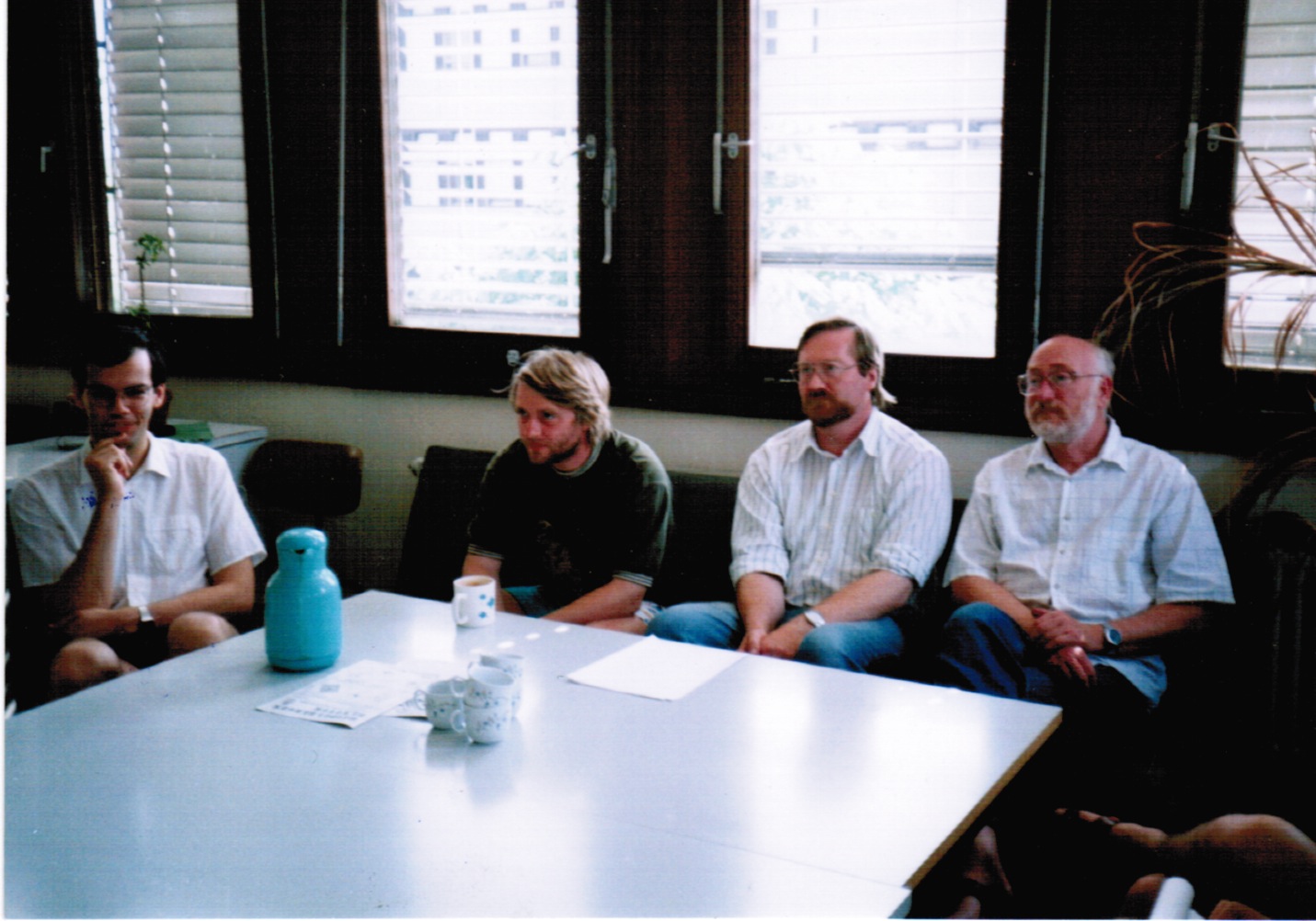}

\quad

\dots and that's what I am still trying to do nowadays. The following is such a discussion on my ignorance in the spirit of Wilhelm's approach to mathematics.

\section{Universal second order constructions}
A non-commutative probability space $(\cA,\ff)$ consists of a unital algebra $\cA$ and a unital linear functional $\ff:\cA\to\CC$. The notion of independence corresponds then to constructions which embed given non-commutative probability spaces into a bigger one; and this construction should be natural or universal in the sense that it does not use any internal structure of the considered non-commutative probability spaces, but works in the same way for all of them. In particular, this means that all elements in an algebra have to be considered equal, there is no special role for generators. In probabilistic terms this corresponds to the requirement that functions of independent variables should also be independent. 

One can formulate the requirements for such constructions in an abstract categorical setting (and this is done by Michael Sch\"urmann \cite{BGS}); but from my more naive approach a universal construction should just be given by formulas to express mixed moments in moments of the individual variables. The universality means then that those formulas cannot use more information than is visible in the mixed moment and, in particular, they must respect associativity. In order to go for the nicest situation I usually, and in particular for today, also restrict to unital situations, where the unit is respected, and to symmetric products, where the product is invariant under exchanging the non-commutative probability spaces. In this setting I could then show \cite{Spe} that the consistency of the calculation rules for mixed moments poses strong conditions on the possibilities and the final result is that there are only two such contructions, namely classical (= tensor) independence and free independence. I think the astonishing point here is not that one can exclude other possibilities, but the fact that this exclusion does not bring you down to one, but to two possibilities. This means of course that at some point your arguments must be of a non-linear nature. Before the birth of free probability \cite{Voi1}, tensor independence was considered as the one and only form of independence in the non-commutative setting. That there is another one came as a surprise and showed that the ideas which one has about what a notion of independence should satisfy are not strong enough to make this unique. Then for a while the pendulum was moving  in the opposite direction and there was the feeling around that if there are two possibilities then there should probably be many more. However, this is not true. There are just two. 

Let us now move forward to second-order probability spaces. There a non-commutative probability space is equipped with an additional object, namely a functional with two arguments
$$\ff_2:\cA\times \cA\to\CC;\qquad (a_1,a_2)\mapsto \ff_2(a_1,a_2),$$
which is linear in each argument and for which the unit also plays a special role, but in the following way:
$$\ff_2(a,1)=0=\ff_2(1,a)$$
for all $a\in\cA$. (Think of $\ff_2$ as a covariance of traces of random matrices.) We also assume that $\ff_2$ is symmetric in its two arguments and that it is tracial in each of its arguments, but I am not sure whether this is relevant for the universality question.
In this setting it is then appropriate to also denote the state $\ff$ from before as $\ff_1$. A second-order non-commutative probability space is thus a triple $(\cA,\ff_1,\ff_2)$. Positivity might be assumed for $\ff_1$ (but is usually not relevant for the arguments), but for $\ff_2$ we don't even know what a good notion of (or replacement for) positivity would be.

Let us now start our game and try to find out what the formulas for mixed moments could look like. We can of course restrict everything to first order, only dealing with $\ff_1$. But then we are back to the old situation and we know that there are only two possibilities. Since we will allow that mixed moments of second order can also depend on moments of first order we have to specify which situation we want to consider for $\ff_1$. Clearly, I vote for the free situation. (Though, it might be interesting to see whether there is a second order theory which goes with tensor independence.)

Now we can look at second order mixed moments. In the following I will assume that I have two or more algebras and I will denote the elements from the first algebra by $a$, the elements from the second algebra by $b$ and so on ...

The first mixed moment of second order is $\ff_2(a,b)$. The only moments in $a$ or in $b$ which I see there are $\ff_1(a)$ and $\ff_1(b)$. Because our formulas must be multi-linear in all appearing arguments, the only possibility is
$$\ff_2(a,b)=\alpha \ff_1(a)\ff_1(b)$$
for some complex number $\alpha\in\CC$. The important point here is of course that $\alpha$ is independent from the concretely considered algebras and elements. So, in particular, we can put
$b=1$, where we get
$$0=\ff_2(a,1)=\alpha \ff_1(a)\ff_1(1)=\alpha\ff_1(a).$$
Since $\ff_1(a)$ can be arbitrary, we
thus have $\alpha=0$. Hence the general formula for this simplest mixed moment must be
$$\ff_2(a,b)=0.$$

Now let's consider a mixed moment of the form $\ff_2(a_1b,a_2)$. The only individual moments I can build from this are $\ff_2(a_1,a_2)\ff_1(b)$, $\ff_1(a_1a_2)\ff_1(b)$, and $\ff_1(a_1)\ff_1(a_2)\ff_1(b)$; hence the formula for the mixed moment must be of the form
$$\ff_2(a_1b,a_2)=\alpha  \ff_2(a_1,a_2)\ff_1(b)+\beta \ff_1(a_1a_2)\ff_1(b)+\gamma \ff_1(a_1)\ff_1(a_2)\ff_1(b)$$
for three coefficients $\alpha,\beta,\gamma\in\CC$.

Putting $b=1$ gives
\begin{align*}
\ff_2(a_1,a_2)&=\alpha  \ff_2(a_1,a_2)\ff_1(1)+\beta \ff_1(a_1a_2)\ff_1(1)+\gamma \ff_1(a_1)\ff_1(a_2)\ff_1(1)\\
&=\alpha  \ff_2(a_1,a_2)+\beta \ff_1(a_1a_2)+\gamma \ff_1(a_1)\ff_1(a_2)
\end{align*}
Since the three different moments in this formula can be choosen independently (here one notices that insisting on positivity is in general not a good idea for such arguments) this
implies that $\alpha=1$ and $\beta=\gamma=0$ and thus 
$$\ff_2(a_1b,a_2)=\ff_2(a_1,a_2)\ff_1(b).$$

Now let's get more serious and consider
\begin{multline*}
\ff_2(a_1b_1,a_2b_2)=\alpha_1\ff_2(a_1,a_2)\ff_2(b_1,b_2)\\+\alpha_2 \ff_2(a_1,a_2)\ff_1(b_1b_2)+\alpha_3\ff_2(a_1,a_2)\ff_1(b_1)\ff_1(b_2)\\
+ \alpha_4\ff_1(a_1a_2)\ff_2(b_1,b_2)+\alpha_5\ff_1(a_1)\ff_1(a_2)\ff_2(b_1,b_2)\\+\alpha_6\ff_1(a_1a_2)\ff_1(b_1b_2)+\alpha_7\ff_1(a_1a_2)\ff_1(b_1)\ff_1(b_2)\\+\alpha_8 \ff_1(a_1)\ff_1(a_2)\ff_1(b_1b_2)+\alpha_9\ff_1(a_1)\ff_1(a_2)\ff_1(b_1)\ff_1(b_2).
\end{multline*}
Putting $b_2=1$ gives
\begin{align*}
\ff_2(a_1,a_2)\ff_1(b_1)&=
\ff_2(a_1b_1,a_2)\\
&=\alpha_2 \ff_2(a_1,a_2)\ff_1(b_1)+\alpha_3\ff_2(a_1,a_2)\ff_1(b_1)\\
&\quad+\alpha_6\ff_1(a_1a_2)\ff_1(b_1)+\alpha_7\ff_1(a_1a_2)\ff_1(b_1)\\
&\quad+\alpha_8 \ff_1(a_1)\ff_1(a_2)\ff_1(b_1)+\alpha_9\ff_1(a_1)\ff_1(a_2)\ff_1(b_1).
\end{align*}
This yields
$$\alpha_2+\alpha_3=1,\qquad \alpha_6+\alpha_7=0,\qquad
\alpha_8+\alpha_9=0.$$
Putting $b_1=1$ gives the same. So let us put now $a_2=1$ instead:
\begin{align*}
\ff_1(a_1)\ff_2(b_1,b_2)&=
\ff_2(a_1b_1,b_2)\\
&=
\alpha_4\ff_1(a_1)\ff_2(b_1,b_2)+\alpha_5\ff_1(a_1)\ff_2(b_1,b_2)\\
&\quad+\alpha_6\ff_1(a_1)\ff_1(b_1b_2)+\alpha_7\ff_1(a_1)\ff_1(b_1)\ff_1(b_2)\\
&\quad+\alpha_8 \ff_1(a_1)\ff_1(b_1b_2)+\alpha_9\ff_1(a_1)\ff_1(b_1)\ff_1(b_2).
\end{align*}
This gives
$$\alpha_4+\alpha_5=1,\qquad \alpha_6+\alpha_8=0,\qquad
\alpha_7+\alpha_9=0.$$
By symmetry in $a$ and $b$ we also know that
$$\alpha_2=\alpha_4,\qquad \alpha_3=\alpha_5,\qquad \alpha_7=\alpha_8.$$
Not all equations are linearly independent, so we do not get a unique solution from them. We need some additional argument(s). For this let us use now associativity of the universal contructions and consider
$\ff_2(a_1b_1c_1,a_2b_2c_2)$ (where three different algebras are now involved). We will see whether and how we can get the term $\ff_1(a_1a_2)\ff_1(b_1b_2)\ff_2(c_1,c_2)$ in the final formula; first we do the calculation according to
\begin{align*}
\ff_2(a_1(b_1c_1),a_2(b_2c_2))&=\alpha_4 \ff_1(a_1a_2) \ff_2(b_1c_1,b_2c_2)+\cdots\\
&=\alpha_4\ff_1(a_1a_2) \alpha_4 \ff_1(b_1b_2)\ff_2(c_1,c_2)+\cdots
\\
&=\alpha_4^2\ff_1(a_1a_2) \ff_1(b_1b_2)\ff_2(c_1,c_2)+\cdots.
\end{align*}
Let us now do the calculation according to
\begin{align*}
\ff_2((a_1b_1)c_1,(a_2b_2)c_2)&=\alpha_4\ff_1(a_1b_1a_2b_2)\ff_2(c_1,c_2)+\cdots\\
&=\alpha_4\cdot 0\cdot \ff_1(a_1a_2)\ff_1(b_1b_2)\ff_2(c_1,c_2)+\cdots.
\end{align*}
In the last line we have used the fact that in first order we have freeness and thus the wanted term does not appear in the factorization of $\ff_1(a_1b_1a_2b_2)$. This calculation thus gives 
$\alpha_4^2=0$ and thus $\alpha_4=0$. From this we then also get
$$\alpha_2=0,\qquad\alpha_3=1, \qquad\alpha_4=0,\qquad \alpha_5=1. $$
Still not enough information, but let us write down what we know up to now:
\begin{multline*}
\ff_2(a_1b_1,a_2b_2)=\alpha\ff_2(a_1,a_2)\ff_2(b_1,b_2)\\+\ff_2(a_1,a_2)\ff_1(b_1)\ff_1(b_2)
+\ff_1(a_1)\ff_1(a_2)\ff_2(b_1,b_2)\\+\beta[\ff_1(a_1a_2)\ff_1(b_1b_2)-\ff_1(a_1a_2)\ff_1(b_1)\ff_1(b_2)\\- \ff_1(a_1)\ff_1(a_2)\ff_1(b_1b_2)+\ff_1(a_1)\ff_1(a_2)\ff_1(b_1)\ff_1(b_2)],
\end{multline*}
where we have renamed $\alpha_1$ and $\alpha_6$ to $\alpha$ and $\beta$, respectively.

We will now follow the term $\ff_2(a_1,a_2)\ff_1(b_1)\ff_1(b_2)\ff_1(c_1)\ff_1(c_2)$ in the same calculation as before. This term can appear on one side as
\begin{align*}
\ff_2((a_1b_1)c_1,(a_2b_2)c_2)&=\ff_2(a_1b_1,a_2b_2)\ff_1(c_1)\ff_1(c_2)+\cdots\\
&=\ff_2(a_1,a_2)\ff_1(b_1)\ff_1(b_2)\ff_1(c_1)\ff_1(c_2)+\cdots
\end{align*}
and on the other side as
\begin{align*}
&\ff_2(a_1(b_1c_1),a_2(b_2c_2))\\&=\alpha\ff_2(a_1,a_2)\ff_2(b_1c_1,b_2c_2)+\ff_2(a_1,a_2)\ff_1(b_1c_1)\ff_1(b_2c_2)+\cdots\\
&=[\alpha\beta+1]\ff_2(a_1,a_2)\ff_1(b_1)\ff_1(b_2)\ff_1(c_1)\ff_1(c_2)+\cdots.
\end{align*}
Comparison of the coefficient gives $1=\alpha\beta+1$, thus $\alpha\beta =0$.
So at least one of $\alpha$ or $\beta$ has to be zero. I am not sure whether there is a meaningful way to normalize the coefficients to 1, but let's do it, and so we essentially are left with two possibilities, namely
\begin{multline*}
\ff_2(a_1b_1,a_2b_2)=\ff_2(a_1,a_2)\ff_1(b_1)\ff_1(b_2)
+\ff_1(a_1)\ff_1(a_2)\ff_2(b_1,b_2)\\+[\ff_1(a_1a_2)\ff_1(b_1b_2)-\ff_1(a_1a_2)\ff_1(b_1)\ff_1(b_2)\\- \ff_1(a_1)\ff_1(a_2)\ff_1(b_1b_2)+\ff_1(a_1)\ff_1(a_2)\ff_1(b_1)\ff_1(b_2)],
\end{multline*}
or 
\begin{multline*}
\ff_2(a_1b_1,a_2b_2)=\ff_2(a_1,a_2)\ff_2(b_1,b_2)\\+\ff_2(a_1,a_2)\ff_1(b_1)\ff_1(b_2)
+\ff_1(a_1)\ff_1(a_2)\ff_2(b_1,b_2).
\end{multline*}
The first case gives us our usual second order freeness. But what about the second case? Is there an argument which excludes the second possibility?
And if we choose the first possibility here, are then the formulas for other mixed moments determined, or do we have to make choices again and again. Experience with the classical case suggests that the concrete form for some mixed moments of small size determines everything, but experience might be misleading. 

This is the point where I am running out of steam and it would be good to sit together with Wilhelm and discuss what those calculations tell us and what to try next, in the hopes of getting some insight into the structure behind all this.

\end{document}